\documentstyle[12pt]{article}
\makeatletter
\renewcommand\thesection{\@Roman\c@section}
\renewcommand\thesubsection{\thesection.\@arabic\c@subsection}
\makeatother

\begin{document}
\begin{titlepage}
\begin{flushright}
math.QA/9905021
\end{flushright}
\vskip.3in

\begin{center}
{\Large \bf Twisted Quantum Affine Superalgebra $U_q[gl(m|n)^{(2)}]$
    and New $U_q[osp(m|n)]$ Invariant R-matrices}
\vskip.3in
{\large Mark D. Gould} and {\large Yao-Zhong Zhang}
\vskip.2in
{\em Department of Mathematics, University of Queensland, Brisbane,
     Qld 4072, Australia

Email: yzz@maths.uq.edu.au}
\end{center}

\vskip 2cm
\begin{center}
{\bf Abstract}
\end{center}

The minimal irreducible representations of $U_q[gl(m|n)]$, i.e. those
irreducible representations that are also irreducible under
$U_q[osp(m|n)]$ are investigated and shown to be affinizable to give
irreducible representations of the twisted quantum affine superalgebra
$U_q[gl(m|n)^{(2)}]$. The $U_q[osp(m|n)]$ invariant R-matrices corresponding 
to the tensor product of any two minimal representations are constructed,
thus extending our twisted tensor product graph method to the 
supersymmetric case. These give new solutions to the spectral-dependent
graded Yang-Baxter equation arising from $U_q[gl(m|n)^{(2)}]$, which exhibit
novel features not previously seen in the untwisted or non-super cases.

\vskip 3cm

\end{titlepage}


\def\a{\alpha}
\def\b{\beta}
\def\d{\delta}
\def\e{\epsilon}
\def\ve{\varepsilon}
\def\g{\gamma}
\def\k{\kappa}
\def\l{\lambda}
\def\o{\omega}
\def\t{\theta}
\def\s{\sigma}
\def\D{\Delta}
\def\L{\Lambda}

\def\G{{\cal G}}
\def\hG{{\hat{\cal G}}}
\def\R{{\cal R}}
\def\hR{{\hat{\cal R}}}
\def\C{{\bf C}}
\def\P{{\bf P}}
\def\Z2{{{\bf Z}_2}}
\def\T{{\cal T}}
\def\H{{\cal H}}
\def\trho{{\tilde{\rho}}}
\def\tphi{{\tilde{\phi}}}
\def\tT{{\tilde{\cal T}}}
\def\uqsnh{{U_q[\widehat{sl(n|n)}]}}
\def\uqs1h{{U_q[\widehat{sl(1|1)}]}}


\def\beq{\begin{equation}}
\def\eeq{\end{equation}}
\def\bea{\begin{eqnarray}}
\def\eea{\end{eqnarray}}
\def\ba{\begin{array}}
\def\ea{\end{array}}
\def\no{\nonumber}
\def\lt{\left}
\def\rt{\right}
\newcommand{\bq}{\begin{quote}}
\newcommand{\eq}{\end{quote}}

\newtheorem{Theorem}{Theorem}
\newtheorem{Definition}{Definition}
\newtheorem{Proposition}{Proposition}
\newtheorem{Lemma}{Lemma}
\newtheorem{Corollary}{Corollary}
\newcommand{\proof}[1]{{\it Proof. }
        #1\begin{flushright}$\Box$\end{flushright}}

\newcommand{\sect}[1]{\setcounter{equation}{0}\section{#1}}
\renewcommand{\theequation}{\thesection.\arabic{equation}}

\sect{Introduction}

The graded Yang-Baxter equation (YBE) 
\beq\label{YB}
R_{12}(z)R_{13}(zw)R_{23}(w)
  =R_{23}(w)R_{13}(zw)R_{12}(z)
\eeq
plays a central role in the theory of supersymmetric quantum integrable systems.
Its solutions $R^{ab}(z)$, usually called R-matrices, 
depend on a spectral parameter 
$z$ and  act on the tensor product
of two graded vector spaces $V_a$ and $V_b$. 
The multiplication in the tensor
product is ${\bf Z}_2$-graded, i.e. for homogeneous elements $x,\;y,\;
x'$ and $y'$,
\beq
(x\otimes y)(x'\otimes y')=(-1)^{[y][x']}\,(x\,x'\otimes y\,y').
\eeq

In one-dimensional lattice integrable systems R-matrices give the 
Hamiltonians of quantum spin chains \cite{Skl79,Kor93}, and 
in statistical mechanics they define the Boltzman weights of
exactly solvable models \cite{Bax82}. 
In integrable quantum field theory R-matrices
give the exact factorizable scattering S-matrices \cite{Zam79,Hol93,Del95c}. 
Therefore, the knowledge of R-matrices has had many physical applications.
This is one of the reasons that the problem of constructing
R-matrices has occupied a fundamental
place in the study of low-dimensional integrable models.

The mathematical framework for the construction of trigonometric
solutions of the (graded) YBE
is given by the quantum affine (super)algebras $U_q({\cal G}^{(k)})$
(see \cite{Jim89} and references therein).
Here for technical reasons we assume $k=1$ or $2$.
Associated to any two finite-dimensional irreducible 
$U_q({\cal G}^{(k)})$-modules $V(\l)$ and $V(\mu)$ there exists a trigonometric
R-matrix $R^{\l,\mu}(z)$ which obeys the (graded) YBE. 
A systematic method, called the tensor product graph (TPG) method, 
for constructing these
R-matrices arising from untwisted quantum affine algebras
was initiated in \cite{Zha91} (see also \cite{Mac92} for the rational cases)
and further developed in \cite{Del94}.
This TPG approach was later generalized
to untwisted quantum affine superalgebras \cite{Del95a}.
In \cite{Del96}, the TPG technique was extended and
a twisted TPG method was developed, which
enables one to construct spectral-dependent R-matrices arising from
twisted quantum affine algebras. The twisted TPG method was used in
\cite{Gan96} to determine soliton S-matrices of certain quantum affine
Toda theories.
By means of the (twisted) TPG method, a large number of
R-matrices have been constructed, leading to many 
new quantum integrable models, quantum spin chains, exactly solvable
lattice models and exact scattering S-matrices. 

It has been an open problem to generalize the TPG method to determine
R-matrices arising from {\it twisted} quantum affine superalgebras. 
In \cite{Gou97}, the $U_q[osp(2|2)]$ invariant R-matrix arising from 
$U_q[gl(2|2)^{(2)}]$ was determined by brute force.  
In this paper we formulate a systematic twisted TPG method, which enables
us to determine R-matrices corresponding to the tensor product of
any two minimal representations of $U_q[gl(m|n)^{(2)}],~(m\leq n,~n>2)$.

The paper is organized as follows: In section II, we review the 
necessary facts about $gl(m|n)$ and its fixed point subalgebra $osp(m|n)$.
In section III, we discuss the twisted quantum affine superalgebra
$U_q[gl(m|n)^{(2)}]$, its minimal representations and give the
equations which uniquely determine the R-matrices (the Jimbo equations).
In sections IV and V we explain how to solve these equations. Our technique is
a supersymmetric generalization of the twisted TPG method introduced
in \cite{Del96}. We have obtained the R-matrices corresponding to
the tensor product $V(\l_a)\otimes V(\l_b)$ of any two minimal
representations $V(\l_a)$ and $V(\l_b)$ of $U_q[gl(m|n)^{(2)}] ~(m\leq n,~
n>2)$. 
Particularly interesting is the case corresponding to
$m=n>2,~a=b$, where an indecomposable representation generally
occurs in the decomposition of
$V(\l_a)\otimes V(\l_a)$, so that the
corresponding R-matrices exhibit new properties not observed in the
untwisted case and in the bosonic case. 
This necessitates an extension of our twisted TPG method.
As we see below 
a combination of the twisted TPG method and the technique used in
\cite{Gou97} is needed to determine the corresponding R-matrices.
Some concluding remarks are given in section VI.

\sect{Priliminaries}

Throughout this paper, 
we assume $n=2r$ is even and set $h=[m/2]$ so that $m=2h$
for even $m$ and $m=2h+1$ for odd $m$.  
For homogeneous operators $A, B$ we use the notation 
$[A, B]=AB-(-1)^{[A][B]}BA$ to denote the usual graded commutator.
Let $E^a_b$ be the standard generators
of $gl(m|n)$ obeying the graded commutation relations
\beq
[E^a_b, E^c_d]=\d^c_b E^a_d-(-1)^{([a]+[b])([c]+[d])}\d^a_d E^c_b.
\eeq
In order to introduce the subalgebra $osp(m|n)$ we first need a graded
symmetric metric tensor $g_{ab}=(-1)^{[a][b]}g_{ba}$ which is assumed 
to be even.
We shall make the convenient choice
\beq
g_{ab}=\xi_a\d_{a\bar{b}},
\eeq
where
\beq
\bar{a}=\left\{
\begin{array}{ll}
m+1-i,~~~& a=i\\
n+1-\mu,~~~& a=\mu,
\end{array}
\right.,~~~~~
\xi_a=\left\{
\begin{array}{ll}
1,~~~& a=i\\
(-1)^\mu,~~~& a=\mu,
\end{array}
\right..
\eeq
In the above equations, $i=1,2,\cdots,m$ and $\mu=1,2,\cdots,n$.
Note that
\beq
\xi_a^2=1,~~~\xi_a\xi_{\bar{a}}=(-1)^{[a]},~~~g^{ab}=\xi_b\d_{a\bar{b}}.
\eeq
As generators of the subalgebra $osp(m|n=2r)$ we take
\beq
\s_{ab}=g_{ac}E^c_b-(-1)^{[a][b]}g_{ac}E^c_a=-(-1)^{[a][b]}\s_{ba}
\eeq
which satisfy the graded commutation relations
\bea
[\s_{ab},\s_{cd}]&=&g_{cb}\s_{ad}-(-1)^{([a]+[b])([c]+[d])}
       g_{ad}\s_{cb}\no\\
& &-(-1)^{[c][d]}\lt(g_{bd}\s_{ac}-(-1)^{([a]+[b])([c]+[d])}
       g_{ac}\s_{db}\rt).
\eea
We have an $osp(m|n)$-module decomposition
\beq
gl(m|n)=osp(m|n)\oplus {\cal T}, ~~~~[{\cal T},{\cal T}]\subset osp(m|n),
\eeq
where ${\cal T}$ is spanned by operators
\beq
T_{ab}=g_{ac} E^c_b+(-1)^{[a][b]}g_{bc} E^c_a=(-1)^{[a][b]}T_{ba}.
\eeq

It is convenient to introduce the Cartan-Weyl generators
\beq
\s^a_b=g^{ac}\s_{cb}=-(-1)^{[a]([a]+[b])}\xi_a\xi_b\s^{\bar{b}}_{\bar{a}}.
\eeq
As a Cartan subalgebra we take the diagonal operators
\beq
\s^a_a=E^a_a-E^{\bar{a}}_{\bar{a}}=-\s^{\bar{a}}_{\bar{a}}.
\eeq
Note that for odd $m=2h+1$ we have $\overline{h+1}=h+1$ and thus
$\s^{h+1}_{h+1}=E^{h+1}_{h+1}-E^{h+1}_{h+1}=0$.

The positive roots of $osp(m|n)$ are given by the even positive roots
(usual positive roots for $o(m)\oplus sp(n)$) together with the odd
positive roots $\d_\mu+\e_i,~1\leq i\leq m,~1\leq\mu\leq r=n/2$,
where we have adopted the useful convention
$\e_{\bar{i}}=-\e_i,~i\leq h=[m/2]$ so that $\e_{h+1}=0$ for odd
$m=2h+1$.  This is consistent with the ${\bf Z}$-gradation
\beq
osp(m|n)=L_{-2}\oplus L_{-1}\oplus L_0\oplus L_1\oplus L_2.
\eeq
Here $L_0=o(m)\oplus gl(r)$, the $gl(r)$ generators are given by
\beq
\s^\mu_\nu=E^\mu_\nu-(-1)^{\mu+\nu}E^{\bar{\nu}}_{\bar{\mu}},~~~
    1\leq \mu,\nu\leq r,
\eeq
and $L_{-2}\oplus L_0\oplus L_2=o(m)\oplus sp(n)$, where $L_2$ gives
rise to an irreducible representation of $L_0$ with highest weight
$(\dot{0}|2,\dot{0})$ spanned by the generators
\beq
\s^\mu_{\bar\nu}=E^\mu_{\bar{\nu}}-\xi_\mu\xi_{\bar{\nu}}E^\nu_{\bar{\mu}}=
    E^\mu_{\bar{\nu}}+(-1)^{\mu+\nu}E^\nu_{\bar{\mu}},~~~1\leq\mu,\nu\leq r.
\eeq
$L_1$ is spanned by odd root space generators
\beq
\s^\mu_i=E^\mu_i+\xi_\mu E^{\bar{i}}_{\bar{\mu}}=E^\mu_i
   +(-1)^\mu E^{\bar{i}}_{\bar{\mu}},~~~1\leq\mu\leq r,~1\leq i\leq m
\eeq
and gives rise to an irreducible representation of $L_0$ with highest
weight $(1,\dot{0}|1,\dot{0})$. $L_{-1},\;L_{-2}$ give rise to
irreducible representations of $L_0$ dual to $L_1,\;L_2$, respectively.
Finally ${\cal T}$ transforms as an irreducible representation of
$osp(m|n)$ [module the first order invariant of $gl(m|n)$] under the
adjoint action with highest weight $\t=\d_1+\d_2$.

The simple roots of $osp(m|n=2r)$ are thus given by the usual (even)
simple roots of $L_0$ together with the odd simple
root $\bar{\a}_s=\bar{\a}_{h+r}=\d_r-\e_1$ 
which is the lowest weight of $L_0$-module $L_1$.
Note that the simple roots of $o(m)$ depend on whether $m$ is odd
or even, and are given here for convenience: For $m=2h$, 
$\bar{\a}_i=\e_i-\e_{i+1}
,~~1\leq i< h,~~~\bar{\a}_h=\e_{h-1}+\e_h$. For $m=2h+1$, 
$\bar{\a}_i=\e_i-\e_{i+1}
,~~1\leq i< h,~~~\bar{\a}_h=\e_h$.  The simple roots of $gl(r)$ are given
by $\bar{\a}_{h+\mu}=\d_\mu-\d_{\mu+1},~~1\leq \mu< r$.
Corresponding to these simple roots we have the generators for $gl(m|n)$ 
(in what follows $s\equiv h+r$):
\bea
&&E_i=\s^i_{i+1},~~~~F_i=\s^{i+1}_i,~~~~H_i=\s^i_i-\s^{i+1}_{i+1},~~~
   1\leq i< h,\no\\
&&E_h=\s^h_{h+1},~~~~F_h=\s^{h+1}_h,~~~~H_h=\s^h_h,~~~{\rm for}~m=2h+1,\no\\
&&E_h=\s^{h-1}_{\overline{h}},~~~~F_h=\s^h_{\overline{h-1}},~~~~
   H_h=\s^{h-1}_{h-1}+\s^h_h,~~~{\rm for}~m=2h,\no\\
&&E_{h+\mu}=\s^\mu_{\mu+1},~~~~F_{h+\mu}=\s^{\mu+1}_\mu,~~~~
   H_{h+\mu}=\s^\mu_\mu-\s^{\mu+1}_{\mu+1},~~~1\leq\mu<r,\no\\
&&E_s=\s^{\mu=r}_{i=1}=E^r_1+(-1)^r E^{\bar{1}}_{\bar{r}},~~~~
  F_s=\s^1_r,~~~~H_s=\s^1_1+\s^r_r,
\eea
which are the standard generators of $osp(m|n)$, together with
\beq
E_0=T^{\mu=\bar{1}}_{\nu=2}=E^{\bar{1}}_2+E^{\bar{2}}_1,~~~~
  F_0=T^2_{\bar{1}},~~~~H_0=-(\s^1_1+\s^2_2),
\eeq
where $E_0$ is the minimal weight vector of ${\cal T}$.
The graded half-sum of the positive roots
of $osp(m|n=2r)$ is given by
\beq
\rho=\frac{1}{2}\sum^h_{i=1}(m-2i)\e_i+\frac{1}{2}\sum^r_{\mu=1}
    (n-m+2-2\mu)\d_\mu.
\eeq

It is convenient to work with generators
\beq
T^a_b=g^{ac}T_{cb}=E^a_b+(-1)^{[a]([a]+[b])}\xi_a\xi_b E^{\bar{b}}
  _{\bar{a}}
\eeq
which should be compared to our $osp(m|n)$ generators $\s^a_b$ above.
This suggests that we introduce the automorphism
\beq
\o(E^a_b)=-(-1)^{[a]([a]+[b])}\xi_a\xi_b E^{\bar{b}}_{\bar{a}}.
\eeq
Then it is easily verified that $\o$ so defined gives an automorphism
of $gl(m|n)$ of order 2, i.e.
\beq
\o\lt([E^a_b, E^c_d]\rt)=[\o(E^a_b),\o(E^c_d)],~~~\o^2=1.
\eeq
Moreover, by inspection, we have
\beq
\o(\s^a_b)=\s^a_b,~~~~\o(T^a_b)=-T^a_b
\eeq
so that $osp(m|n=2r)$ is the fixed point subalgebra of $\o$ while ${\cal T}$
corresponds to the eigenspace corresponding to eigenvalue $-1$ of $\o$.

\sect{Twisted quantum affine superalgebra $U_q[gl(m|n)^{(2)}]$}

Following the usual notation, we now set
\beq
\hat{L}\equiv gl(m|n),~~~\hat{L}_0\equiv osp(m|n),~~~
\hat{L}_1\equiv {\cal T}.
\eeq
With this notation we have
\beq
\hat{L}=\hat{L}_0\oplus \hat{L}_1.
\eeq
The corresponding twisted affine superalgebra $\hat{L}^{(2)}$ admits
the decomposition
\beq
\hat{L}^{(2)}=\bigoplus_{l\in {\bf Z}/2} \hat{L}^{(2)}_l\oplus
  {\bf C}c\oplus{\bf C} d,
\eeq
where
\beq
\hat{L}^{(2)}_m=\lt\{
\begin{array}{ll}
\hat{L}_0(l), & ~~~l\in {\bf Z}\\
\hat{L}_1(l), & ~~~l\in {\bf Z}+\frac{1}{2},
\end{array}
\rt.
\eeq
with
$\hat{L}_t(l)=\{x(l)~|~x\in\hat{L}_t\},~  t=0,1$.
The graded commutation relations are defined by
\bea
&&[x(l),y(l')]=[x,y](l+l')+l\, c\,\d_{l+l',0}(x,y),\no\\
&&[d,x(l)]=l x(l),~~~[c,x(l)]=[c,d]=0,
\eea
where $(~,~)$ is a fixed invariant bi-linear form on $\hat{L}$.
As a Cartan subalgebra of
$\hat{L}^{(2)}$, we take
\beq
\hat{\cal H}=H(0)\oplus {\bf C}c\oplus{\bf C}d,
\eeq
where $H\in\hat{L}_0$ is a Cartan subalgebra of $\hat{L}_0$. As a set of
simple roots of $\hat{L}^{(2)}$, we take
\beq
\a_0=-\t+\frac{1}{2}\d,~~~\a_i=\bar{\a}_i,~~1\leq i\leq s,
\eeq
where $\t=\d_1+\d_2$, $\bar{\a}_i$ are the simple roots of $\hat{L}_0$
and $\frac{1}{2}\d$ is the minimal positive imaginary root.

Then $\hat{L}^{(2)}$ admits a standard set of simple generators $\{e_i,~
f_i,~h_i| 0\leq i\leq s\}$, given by
\bea
&&e_i=E_i(0),~~~~f_i=F_i(0),~~~~h_i=H_i(0),~~~1\leq i\leq s,\no\\
&&e_0=E_0(\frac{1}{2}),~~~~f_0=F_0(-\frac{1}{2}),~~~~h_0=H_0(0)+\frac{1}{2}c.
\eea
They satisfy the following relations
\bea
&&[h_i,e_j]=(\a_i,\a_j)e_j,~~~~
   [h_i,f_j]=-(\a_i,\a_j)f_j,\no\\
&& {[e_i, f_j]}=\d_{ij}h_i,\no\\
&&\lt({\rm ad}e_i\rt)^{1-a_{ij}}e_j=0=
  \lt({\rm ad}f_i\rt)^{1-a_{ij}}f_j,~~~i\neq j,
\eea
where $a_{ij}$ is the Cartan matrix
\beq
a_{ij}=\lt\{
\begin{array}{ll}
\frac{2(\a_i,\a_j)}{(\a_i,\a_i)}, & ~~~i\neq s\\
(\a_i,\a_j), & ~~~i=s.
\end{array}
\rt.
\eeq
We remark that there are also higher order Serre relations \cite{Yam96}.
Since these extra Serre relations are not important for 
our purpose here, we shall ignore them in this paper.

Twisted quantum affine superalgebra $U_q(\hat{L}^{(2)})$ is a
$q$-deformation of the universal enveloping algebra 
$U(\hat{L}^{(2)})$ of $\hat{L}^{(2)}$. In terms of the simple
generators $\{e_i,~f_i,~h_i | 0\leq i\leq s\}$, the
defining relations of $U_q(\hat{L}^{(2)})$ are given by
\bea
&&[h_i,e_j]=(\a_i,\a_j)e_j,~~~~
   [h_i,f_j]=-(\a_i,\a_j)f_j,\no\\
&& {[e_i, f_j]}=\d_{ij}\frac{q^{h_i}-q^{-h_i}}{q-q^{-1}},\no\\
&&\lt({\rm ad}_{q_i}e_i\rt)^{1-a_{ij}}e_j=0=
  \lt({\rm ad}_{q_i}f_i\rt)^{1-a_{ij}}f_j,~~~i\neq j,
\eea
where $q_i=q^{(\a_i,\a_i)/2}$. Here again we have ignored the extra
$q$-Serre relations.

$U_q(\hat{L}^{(2)})$ is a quasi-triangular Hopf superalgebra with coproduct
$\D$ given by
\bea
&&\D(q^{\pm h_i/2})=q^{\pm h_i/2}\otimes q^{\pm h_i/2},\no\\
&&\D(e_i)=q^{-h_i/2}\otimes e_i+ e_i\otimes q^{h_i/2},\no\\
&&\D(f_i)=q^{-h_i/2}\otimes f_i+ f_i\otimes q^{h_i/2}.
\eea
Throughout ${\cal R}$ denotes the universal R-matrix of
$U_q(\hat{L}^{(2)})$ which by definition satisfies
\bea
&&
{\cal R}\D(a)=\D^T(a){\cal R},~~\forall a\in U_q(\hat{L}^{(2)}),
\no\\
 &&
(1\otimes\D){\cal R}={\cal R}_{13}{\cal R}_{12},~~~
(\D\otimes1){\cal R}={\cal R}_{13}{\cal R}_{23}\label{rdef}
\eea
where $\D^T(a)$ is the opposite coproduct. A direct consequence
of the above relations is that ${\cal R}$ satisfies the graded YBE.
\beq
{\cal R}_{12}{\cal R}_{13}{\cal R}_{23}=
{\cal R}_{23}{\cal R}_{13}{\cal R}_{12}
\eeq
Note that the generators $q^{\pm h_i/2}, e_i, f_i,~(1\leq i\leq s)$
generate the quantum algebra $U_q(\hat{L}_0)$ which is a
quasi-triangular Hopf subsuperalgebra of $U_q(\hat{L}^{(2)})$. We denote 
by $R$ the universal R-matrix of $U_q(\hat{L}_0)$.

We now consider finite-dimensional irreducible
representations of $U_q(\hat{L}^{(2)})$
in order to give new solutions to the graded YBE.
If a finite-dimensional irreducible representation 
of $U_q(\hat{L}^{(2)})$ can be constructed by affinizing an
irreducible  representation
of $U_q(\hat{L}_0)$, we say this irreducible
representation of $U_q(\hat{L}_0)$ is
affinizable. Not all irreducible
representations of $U_q(\hat{L}_0)$ are affinizable
to give representations of $U_q(\hat{L}^{(2)})$. An important role  is
played by minimal irreducible representations of $\hat{L}$ which are
those irreducible representations which are also irreducible under
$\hat{L}_0$. We shall see below
that all such irreducible representations of $\hat{L}_0$
can be quantized to give irreducible representations of $U_q(\hat{L}_0)$
which are affinizable to
$U_q(\hat{L}^{(2)})$.  Using the graded fermion formalism developed in
our previous paper \cite{Gou99a}, one can show that
the anti-symmetric tensor irreducible
representations of $\hat{L}$, with highest weights
\beq
\L_a=\lt\{
\begin{array}{ll}
(\dot{1}_a, \dot{0}|\dot{0}), &~~~a\leq m,\\
(\dot{1}|a-m,\dot{0}), & ~~~a>m \label{min-weight}
\end{array}
\rt.
\eeq
are minimal, and thus give rise to irreducible modules under
$U_q(\hat{L}_0)$, with 
the corresponding $U_q(\hat{L}_0)$ highest weight
\beq
\l_a=(\dot{0}|a,\dot{0}).
\eeq
We thus denote this minimal irreducible
representation by $V(\l_a)$ to denote its highest weight under
$U_q(\hat{L}_0)$.

We define an automorphism $D_z$ of $U_q(\hat{L}^{(2)})$ by
\beq
D_z(e_i)=z^{\d_{i0}}e_i,~~~~
D_z(f_i)=z^{-\d_{i0}}f_i,~~~~
D_z(h_i)=h_i.
\eeq
Given any two minimal irreps $\pi_\l$ and $\pi_\mu$ of $U_q(\hat{L}_0)$
and their affinizations to irreducible representations of $U_q(\hat{L}^{(2)})$,
we obtain a one-parameter family of representations 
$\D_{\l\mu}^z$ of $U_q(\hat{L}^{(2)})$ on
$V(\l)\otimes V(\mu)$ defined by
\beq
\D_{\l\mu}^z(a)=\pi_\l\otimes \pi_\mu\left((D_z\otimes 1)
\D(a)\right),~~~\forall a\in U_q(\hat{L}^{(2)}),
\eeq
where $z$ is the spectral parameter. We define the spectral
parameter dependent R-matrix
\beq
R^{\l\mu}(z)=(\pi_\l\otimes\pi_\mu)\left((D_z\otimes 1){\cal R}\right).
\eeq
It follows that this R-matrix gives a solution to
the spectral parameter dependent YBE.
{}From the defining property (\ref{rdef}) of the universal R-matrix
one derives the equations
\beq\label{Jimbo}
R^{\l\mu}(z)\,\D_{\l\mu}^z(a)=(\D^T)_{\l\mu}^z(a)\,R^{\l\mu}(z)
\eeq
which, because the representations $\D_{\l\mu}^z$
are irreducible for generic $z$, 
uniquely determine $R^{\l\mu}(z)$ up to a scalar function 
of $z$. 

We normalize $R^{\l\mu}(u)$ such that
\beq
\check{R}^{\l\mu}(z)\check{R}^{\mu\l}(z^{-1})=I~~
\mbox{and}~~R(0)=\pi_\l\otimes\pi_\mu(R),
\eeq
where $R$ is the R-matrix of $U_q(\hat{L}_0)$ and 
$\check{R}^{\l\mu}(z)=P\,R^{\l\mu}(z)$ with
$P:V(\l)\otimes V(\mu)\rightarrow V(\mu)\otimes V(\l)$
the usual graded permutation operator.

In order for the equation (\ref{Jimbo}) to hold for 
all $a\in U_q(\hat{L}^{(2)})$ it is sufficient that it holds for all
$a\in U_q(\hat{L}_0)$ and in addition for the extra generator
$e_0$.
The relation for $e_0$ reads explicitly
\bea
&&
R^{\l\mu}(z)
\left(z\,\pi_\l(e_0)\otimes\pi_\mu(q^{h_0/2})+
\pi_\l(q^{-h_0/2})\otimes \pi_\mu(e_0)\right)
\no\\
&&
~~~~~=\left(z\pi_\l(e_0)\otimes\pi_\mu(q^{h_0/2})+
\pi_\l(q^{-h_0/2})\otimes\pi_\mu(e_0)\right)
R^{\l\mu}(z),
\eea
or equivalently
\bea
&&
\check{R}^{\l\mu}(z)
\left(z\,\pi_\l(e_0)\otimes\pi_\mu(q^{h_0/2})+
\pi_\l(q^{-h_0/2})\otimes \pi_\mu(e_0)\right)
\no\\
&&
~~~~~=\left(\pi_\mu(e_0)\otimes\pi_\l(q^{h_0/2})+
z\,\pi_\mu(q^{-h_0/2})\otimes\pi_\l(e_0)\right)
\check{R}^{\l\mu}(z).\label{jimbo-eq}
\eea
Eq.(\ref{jimbo-eq}) is the Jimbo equation for the twisted quantum
affine superalgebras. In the next section,
we shall solve this equation for all minimal representations
of $U_q(\hat{L}^{(2)})$, for the case $m\leq n,~n>2$.

\sect{Twisted Tensor Product Graph}

We shall determine $U_q[osp(m|n)]$ invariant R-matrices 
$\check{R}^{\l_a\l_b}(z)$ on the tensor product of any two minimal
irreducible representations $V(\l_a)\otimes V(\l_b),~a\leq b$, arising from
$U_q[gl(m|n)^{(2)}]$. 

Let $V(\l)$ and $V(\mu)$ denote any two minimal irreducible representations of
$U_q(\hat{L}^{(2)})$. Except for the case $a-b=m-n=0$ (see below), the tensor product
module $V(\l)\otimes V(\mu)$ is completely reducible into
irreducible $U_q(\hat{L}_0)$-modules as
\beq\label{dec}
V(\l)\otimes V(\mu)=\bigoplus_\nu V(\nu)
\eeq
and there are no multiplicities in this decomposition.
We denote by $P_\nu^{\l\mu}$  the projection operator of
$V(\l)\otimes V(\mu)$ onto $V(\nu)$ and set
\beq
{\bf P}^{\l\mu}_\nu=\check{R}^{\l\mu}(1)\,P^{\l\mu}_\nu=
P^{\mu\l}_\nu\,\check{R}^{\l\mu}(1).
\eeq
We may thus write
\beq
\check{R}^{\l\mu}(z)=\sum_\nu\,\rho_\nu(z)\,
  {\bf P}^{\l\mu}_\nu,~~~
  \rho_\nu(1)=1.
\eeq
Following our previous approach \cite{Del94}, the coefficients 
$\rho_\nu(z)$ may be determined according to the recursion
relation
\beq\label{rec}
\rho_\nu(z)=\frac{q^{C(\nu)/2}+\e_\nu\e_{\nu'}z\,q^{C(\nu')/2}}
{z\,q^{C(\nu)/2}+\e_\nu\e_{\nu'}\,q^{C(\nu')/2}}
\rho_{\nu'}(z),
\eeq
which holds for any $\nu\neq\nu'$ for which
\beq\label{edge}
P^{\l\mu}_\nu\left(\pi_\l(e_0)\otimes\pi_\mu(q^{h_0/2})\right)
P^{\l\mu}_{\nu'}\neq 0.
\eeq
Here $C(\nu)$ is the eigenvalue of the universal Casimir element
of $\hat{L}_0$ on $V(\nu)$ and $\e_\nu$ denotes the parity
of $V(\nu)\subseteq V(\l)\otimes V(\mu)$.

To graphically encode the recursive relations between the 
different $\rho_\nu$ we introduce the {\bf Twisted TPG}
$\tilde{G}^{\l\mu}$ associated to the tensor
product module $V(\l)\otimes V(\mu)$. The nodes of this
graph are given by the highest weights $\nu$ of the
$U_q(\hat{L}_0)$-modules occuring in the decomposition (\ref{dec})
of the tensor product module. There is an edge between two
nodes $\nu\neq\nu'$ iff (\ref{edge}) holds.

Given a tensor product module and its decomposition, it is not 
in general an easy task to determine the twisted TPG
 because in order to determine between which nodes of the
graph relation (\ref{edge}) holds requires detailed calculations. 
We therefore introduce the
{\bf Extended Twisted TPG} $\tilde{\G}^{\l\mu}$
which has the same set of nodes as the twisted TPG
but has an edge between two vertices $\nu\neq\nu'$ whenever
\beq\label{edge1}
V(\nu')\subseteq V(\theta)\otimes V(\nu)
\eeq
and
\beq\label{edge2}
\epsilon_\nu \epsilon_{\nu'}=
\left\{\begin{array}{l}
+1~~~\mbox{if }V(\nu)\mbox{ and }V(\nu')\mbox{ are in the same
irreducible representation of }\hat{L}\\
-1~~~\mbox{if }V(\nu)\mbox{ and }V(\nu')\mbox{ are in different
irreducible representations of }\hat{L}.\end{array}\right.
\eeq
The conditions (\ref{edge1}) and (\ref{edge2}) are necessary
conditions for (\ref{edge}) to hold and therefore the twisted TPG
is contained in the extended twisted TPG. 
To see why (\ref{edge1}) is a necessary condition for
(\ref{edge}) one must
realize that $e_0\otimes q^{h_0/2}$ is the lowest component of
a tensor operator corresponding to $V(\theta)$.
The necessity of (\ref{edge2}) follows 
from the fact 
that two vertices $\nu\neq\nu'$ connected by an edge in the
twisted TPG (i.e., for which (\ref{edge}) is satisfied)
must have the {\em same} parity if $V(\nu)$ and $V(\nu')$
belong to the {\em same} irreducible $\hat{L}$-module while they
must have {\em opposite} parities if they belong to different
irreducible $\hat{L}$-modules.

While the extended twisted TPG will always include
the twisted TPG, it will in general have more
edges. Only if the extended twisted TPG is a tree are we guaranteed that it 
coincides with the twisted TPG.

We will impose a relation (\ref{rec}) for every edge in the
extended twisted TPG. Because the extended
TPG will in general have more edges than the
(unextended) twisted TPG, we will be imposing
too many relations. These relations may be inconsistent and we
are therefore not guaranteed a solution.
If however a solution exists, then it must be the
unique correct solution to the Jimbo's equation.

As seen below, for the cases we are considering, the 
extended twisted TPG is always consistent and
thus will always give rise to a solution of the graded YBE.

Throughout we adopt the convenient notation
\beq
<a>_\pm=\frac{1\pm z\,q^a}{z\pm q^a},
\eeq
so that the relation (\ref{rec}) may be expressed as
\beq
\rho_\nu(z)=\left\langle\frac{C(\nu')-C(\nu)}{2}\right\rangle
_{\epsilon_\nu\epsilon_{\nu'}}\,
\rho_{\nu'}(z).
\eeq

\sect{Solution to Jimbo Equation}

We will now determine the R-matrices for any tensor product
$V(\l_a)\otimes V(\l_b)$ of two minimal representations 
$V(\l_a),~V(\l_b)$ of $U_q(\hat{L}^{(2)})$.  Note that
for our case, $\t=\d_1+\d_2$.

\subsection{The case of either $a\leq b,~m<n,~n>2$ or $a<b,~m=n>2$}

We consider the decomposition of the tensor product of two minimal
irreducible representations of $U_q(\hat{L}_0)$: $V(\l_a)\otimes V(\l_b)$.
Recall that $V(\l_a),~V(\l_b)$  also carry irreducible representations
of $\hat{L}$. 
We denote by $\hat{V}(a,b)$ the ``two-column" irreducible representation
of  $\hat{L}$ with highest weight
\beq
{\L}_{a,b}=\lt\{
\begin{array}{ll}
(\dot{2}_a,\dot{1}_b,\dot{0}|\dot{0}),~~~& a+b\leq m\\
(\dot{2}_a,\dot{1}|a+b-m,\dot{0}),~~~& a\leq m,~a+b>m\\
(\dot{2}|a+b-m,a-m,\dot{0}),~~~& a>m.
\end{array}
\rt. 
\eeq
Then we have the following decomposition
 into irreducible $\hat{L}$-modules
\bea
V(\l_a)\otimes V(\l_b)&\equiv& \hat{V}({\L}_a)
  \otimes \hat{V}({\L}_b)\no\\
&=&\bigoplus^a_{c=0}\hat{V}(c, a+b-2c).
\eea

The reduction of the two-column irreducible representations of $\hat{L}$
into irreducible representations of $\hat{L}_0$ has been worked out
in our previous paper \cite{Gou99a} by using the quasi-spin graded-fermion
formalism. We thus arrive at the following irreducible $\hat{L}_0$
, and thus $U_q(\hat{L}_0)$ module decomposition,
\beq
V(\l_a)\otimes V(\l_b)
   =\bigoplus^a_{c=0}\bigoplus^c_{k=0}V(k, a+b-2c);\label{vv-decom}
\eeq
here and throughout $V(a,b)$ denotes an irreducible $U_q(\hat{L}_0)$
module with highest weight
\beq
\l_{a,b}=(\dot{0}|a+b,a,\dot{0}).
\eeq
Note that one can only get an indecomposable in (\ref{vv-decom}) when
$m=n>2$ and $a+b-2c=0$, which is the case we shall consider in
the next subsection. Since $a\leq b,~c\leq a$, this can only occur when
$a=b$ and $c=a$. In that case the $U_q(\hat{L}_0)$-modules $V(k,0),~
k=0,1,$ will form an indecomposable. 

We now show that the minimal irreducible $U_q(\hat{L}_0)$,
$V(\l_a)$, 
with highest weight $\l_a$, are affinizable. We first
note that the vector module $V(\l_1)$ of $U_q(\hat{L}_0)$ with higest
weight $\l_1=(\dot{0}|1,\dot{0})$ is minimal. It is also affinizable 
since it is undeformed.  

Following our previous approach \cite{Del96}, 
we consider the corresponding {\it twisted} TPG for
$V(\l_1)\otimes V(\l_1)$
\bea\label{ttpg1}
\unitlength=1mm
\linethickness{0.4pt}
\begin{picture}(102.60,12.60)(40,13)
\put(60.00,20.00){\circle*{5.20}}
\put(80.00,20.00){\circle*{5.20}}
\put(60.00,16.00){\makebox(0,0)[ct]{$\l_0$}}
\put(80.00,16.00){\makebox(0,0)[ct]{$\l_1+\d_2$}}
\put(100.00,20.00){\circle*{5.20}}
\put(100.00,16.00){\makebox(0,0)[ct]{$\l_2$}}
\put(60.00,20.00){\line(1,0){40.00}}
\put(60.00,24.00){\makebox(0,0)[cb]{$-$}}
\put(80.00,24.00){\makebox(0,0)[cb]{$-$}}
\put(100.00,24.00){\makebox(0,0)[cb]{$+$}}
\end{picture}
\eea
where $\pm$ indicate the parities. 
Since $\l_2$ is an extremal node on the twisted TPG it
follows that $V(\l_2)$ is affinizable \cite{Del95b}, i.e. it too
carries an irreducible representation of $U_q(\hat{L}^{(2)})$. 
Note that the twisted TPG (\ref{ttpg1})
has a quite different topology to the untwisted TPG:
\bea\no
\unitlength=1mm
\linethickness{0.4pt}
\begin{picture}(102.60,12.60)(40,13)
\put(60.00,20.00){\circle*{5.20}}
\put(80.00,20.00){\circle*{5.20}}
\put(60.00,16.00){\makebox(0,0)[ct]{$\l_0$}}
\put(80.00,16.00){\makebox(0,0)[ct]{$\l_2$}}
\put(100.00,20.00){\circle*{5.20}}
\put(100.00,16.00){\makebox(0,0)[ct]{$\l_1+\d_2$}}
\put(60.00,20.00){\line(1,0){40.00}}
\put(60.00,24.00){\makebox(0,0)[cb]{$-$}}
\put(80.00,24.00){\makebox(0,0)[cb]{$+$}}
\put(100.00,24.00){\makebox(0,0)[cb]{$-$}}
\end{picture}
\eea
which indicates that $V(\l_2)$ is not affinizable to carry an irreducible
representation of $U_q(\hat{L}_0^{(1)})$.

More generally, we have the following twisted TPG for for 
$V(\l_1)\otimes V(\l_a)$
\bea\label{ttpg2}
\unitlength=1mm
\linethickness{0.4pt}
\begin{picture}(102.60,12.60)(40,13)
\put(60.00,20.00){\circle*{5.20}}
\put(80.00,20.00){\circle*{5.20}}
\put(100.00,20.00){\circle*{5.20}}
\put(60.00,20.00){\line(1,0){40.00}}
\put(60.00,16.00){\makebox(0,0)[ct]{$\l_{a-1}$}}
\put(80.00,16.00){\makebox(0,0)[ct]{$\l_a+\d_2$}}
\put(100.00,16.00){\makebox(0,0)[ct]{$\l_{a+1}$}}
\put(60.00,24.00){\makebox(0,0)[cb]{$-$}}
\put(80.00,24.00){\makebox(0,0)[cb]{$-$}}
\put(100.00,24.00){\makebox(0,0)[cb]{$+$}}
\end{picture}
\eea
so that $\l_{a+1}$ is an extremal node and hence, by recursion,
each of the irreducible representations
$V(\l_a)$, is affinizable to $U_q(\hat{L}^{(2)})$. 
Thus we have shown that
all minimal irreducible $U_q(\hat{L}_0)$ modules $V(\l_a)$ are affinizable to carry
irreducible representations of $U_q(\hat{L}^{(2)})$.

The $U_q(\hat{L}_0)$-module decomposition of the tensor product of any
two such representations is given by (\ref{vv-decom}). Thus we have
the corresponding extended twisted TPG for $V(\l_a)\otimes V(\l_b)$,
given below by Figure 1.

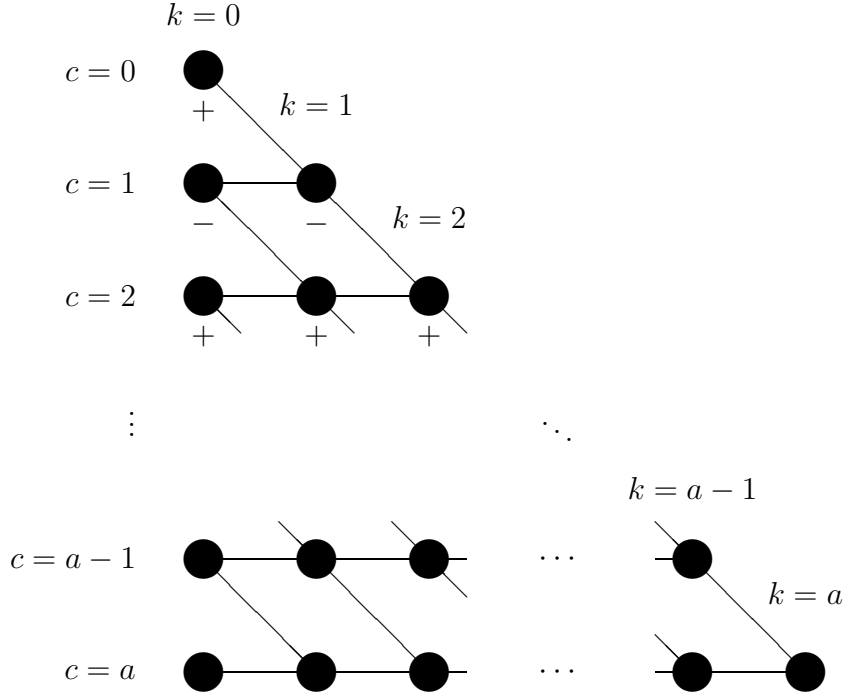
\begin{figure}[ht]
\unitlength=1.00mm
\linethickness{0.4pt}
\begin{picture}(129.60,92.00)
\put(47.00,86.00){\line(1,-1){35.00}}
\put(77.00,56.00){\line(-1,0){30.00}}
\put(67.00,51.00){\line(-1,1){20.00}}
\put(47.00,71.00){\line(1,0){15.00}}
\put(47.00,56.00){\line(1,-1){5.00}}
\put(47.00,6.00){\line(1,0){35.00}}
\put(47.00,21.00){\line(1,0){35.00}}
\put(47.00,21.00){\line(1,-1){15.00}}
\put(77.00,6.00){\line(-1,1){20.00}}
\put(127.00,6.00){\line(-1,0){20.00}}
\put(112.00,6.00){\line(-1,1){5.00}}
\put(127.00,6.00){\line(-1,1){20.00}}
\put(112.00,21.00){\line(-1,0){5.00}}
\put(72.00,26.00){\line(1,-1){10.00}}
\put(47.00,86.00){\circle*{5.20}}
\put(47.00,71.00){\circle*{5.20}}
\put(62.00,71.00){\circle*{5.20}}
\put(62.00,56.00){\circle*{5.20}}
\put(47.00,56.00){\circle*{5.20}}
\put(77.00,56.00){\circle*{5.20}}
\put(77.00,21.00){\circle*{5.20}}
\put(62.00,21.00){\circle*{5.20}}
\put(47.00,21.00){\circle*{5.20}}
\put(47.00,6.00){\circle*{5.20}}
\put(62.00,6.00){\circle*{5.20}}
\put(77.00,6.00){\circle*{5.20}}
\put(112.00,6.00){\circle*{5.20}}
\put(127.00,6.00){\circle*{5.20}}
\put(112.00,21.00){\circle*{5.20}}
\put(47.00,82.00){\makebox(0,0)[ct]{$+$}}
\put(47.00,67.00){\makebox(0,0)[ct]{$-$}}
\put(47.00,52.00){\makebox(0,0)[ct]{$+$}}
\put(62.00,67.00){\makebox(0,0)[ct]{$-$}}
\put(62.00,52.00){\makebox(0,0)[ct]{$+$}}
\put(77.00,52.00){\makebox(0,0)[ct]{$+$}}
\put(94.00,21.00){\makebox(0,0)[cc]{$\cdots$}}
\put(94.00,6.00){\makebox(0,0)[cc]{$\cdots$}}
\put(38.00,86.00){\makebox(0,0)[rc]{$c=0$}}
\put(38.00,71.00){\makebox(0,0)[rc]{$c=1$}}
\put(38.00,56.00){\makebox(0,0)[rc]{$c=2$}}
\put(38.00,21.00){\makebox(0,0)[rc]{$c=a-1$}}
\put(38.00,6.00){\makebox(0,0)[rc]{$c=a$}}
\put(62.00,80.00){\makebox(0,0)[cb]{$k=1$}}
\put(77.00,65.00){\makebox(0,0)[cb]{$k=2$}}
\put(112.00,29.00){\makebox(0,0)[cb]{$k=a-1$}}
\put(127.00,15.00){\makebox(0,0)[cb]{$k=a$}}
\put(38.00,40.00){\makebox(0,0)[rc]{$\vdots$}}
\put(94.00,39.00){\makebox(0,0)[cc]{$\ddots$}}
\put(47.00,92.00){\makebox(0,0)[cb]{$k=0$}}
\end{picture}
\caption{\small The extended twisted TPG for $U_q[gl(m|n)^{(2)}]$
for the tensor product $V(\l_a)
\otimes V(\l_b)$, for the cases we are considering. 
The vertex labelled by the pair $(c,k)$ corresponds to the
irreducible $U_q(\hat{L}_0)$ module $V(k, a+b-2c)$. The rows and columns are
labelled by $c$ and $k$, respectively. The rows labelled
by $c$ correspond to irreducible representations of $\hat{L}_0$ in
the same irreducible representation of $\hat{L}$ and so all vertices
in a row have the same parity. 
The vertices along the diagonal edges of
the graph have alternating parities.
\label{fig2}}
\end{figure}

To see that this graph is consistent we have to consider the
closed loops of the form
\bea
\unitlength=1.30mm
\linethickness{0.4pt}
\begin{picture}(93.00,28.00)(30,20)
\put(60.00,40.00){\circle*{4.00}}
\put(75.00,40.00){\circle*{4.00}}
\put(75.00,25.00){\circle*{4.00}}
\put(89.00,25.00){\circle*{4.00}}
\put(60.00,40.00){\line(1,0){15.00}}
\put(75.00,40.00){\line(1,-1){14.00}}
\put(89.00,26.00){\line(-1,0){14.00}}
\put(75.00,26.00){\line(-1,1){15.00}}
\put(56.00,40.00){\makebox(0,0)[rc]{$(c,k-1)$}}
\put(79.00,40.00){\makebox(0,0)[lc]{$(c,k)$}}
\put(93.00,25.00){\makebox(0,0)[lc]{$(c+1, k+1)$}}
\put(71.00,25.00){\makebox(0,0)[rc]{$(c+1,k)$}}
\put(60.00,37.00){\makebox(0,0)[ct]{$+$}}
\put(75.00,37.00){\makebox(0,0)[ct]{$+$}}
\put(75.00,22.00){\makebox(0,0)[ct]{$-$}}
\put(89.00,22.00){\makebox(0,0)[ct]{$-$}}
\end{picture}\no
\eea
where we have indicated the relative parities of the vertices.
We denote the eigenvalue of the universal Casimir element of $\hat{L}_0$
on the irreducible representation labelled by $(c,k)$ by $C_{c,k}$. Then
it is easily seen that
\bea
&&
C_{c,k}-C_{c,k-1}=C_{c+1,k+1}-C_{c+1,k}=2(\rho,\d_1+\d_2)-2(a+b-1)+4(c-k),
\no\\
&&
C_{c,k}-C_{c+1,k+1}=C_{c,k-1}-C_{c+1,k}=2(\rho,\d_1-\d_2)-2(a+b-1)+2c.
\eea
This implies that the extended twisted TPG is consistent, i.e. that
the recursion relations (\ref{rec}) give the same result independent of the
path along which one recurses.

We can now read off the R-matrix from the extended twisted TPG 
\bea
\check{R}^{\l_a,\l_b}(z)&=&\sum_{c=0}^a\sum_{k=0}^c\prod^{c-k}_{j=1}
  \langle m-n+2j-a-b\rangle_+\no\\
& & \prod^c_{i=1}
  \langle i-a-b-1 \rangle_-\,{\bf P}^{\l_a,\l_b}_{(a+b-2c+k)\d_1+k\d_2}
\eea

\subsection{The case of $a=b,~m=n>2$}

As indicated in the last subsection, for the case at hand the
$U_q[osp(n|n)]$-modules $V(k,0)$, $~k=0,1$, appearing in the right hand
side of the tensor product decomposition (\ref{vv-decom}), form an
indecomposable representation of $U_q[osp(n|n)]$. From now on we
denote by $V$ this indecomposable module. We thus have an
$U_q[osp(n|n)]$ module decomposition
\beq
V(\l_a)\otimes V(\l_a)=\bigoplus_{\nu} V(\nu)\bigoplus V,\label{decom-indecom}
\eeq
where the sum on $\nu$ is over the irreducible highest weights and
$V$ is the unique indecomposable.  Note that $V$ contains a unique
submodule $\bar{V}(\d_1+\d_2)$ which is maximal, indecomposable and
cyclically generated by a maximal vector of weight $\d_1+\d_2$ such
that $V/\bar{V}(\d_1+\d_2)\cong V(\dot{0}|\dot{0})$ (the trivial
$U_q[osp(n|n)]$-module). Moreover $V$ contains a unique irreducible
submodule $V(\dot{0}|\dot{0})\subset \bar{V}(\d_1+\d_2)$. The usual form
of Schur's lemma applies to $\bar{V}(\d_1+\d_2)$  and so the space
of $U_q[osp(n|n)]$ invariants in ${\rm End}(V)$ has dimension 2 (see
Appendix A). It is
spanned by the identity operator $I$ together with an invariant $N$
(unique up to scalar multiples) satisfying
\beq
N\,V=V(\dot{0}|\dot{0})\subset\bar{V}(\d_1+\d_2),~~~~
N\,\bar{V}(\d_1+\d_2)=(0).\label{N-action}
\eeq
It follows that $N$ is nilpotent, i.e.
\beq
N^2=0.
\eeq

We now determine the extended twisted TPG for the decomposition given
by (\ref{decom-indecom}). We note that  $V$ can only be connected to two 
nodes corresponding to highest weights
\beq
\nu=\lt\{
\begin{array}{ll}
2\d_1~~({\rm opposite~ parity}), &~~~(c,k)=(a-1,0)\\
2(\d_1+\d_2)~~({\rm same~ parity}), &~~~(c,k)=(a,2).\label{nu-node}
\end{array}
\rt.
\eeq
We thus arrive at the (consistent)
extended twisted TPG for (\ref{decom-indecom}), given by Figure 2.

\begin{figure}[ht]
\unitlength=1.00mm
\linethickness{0.4pt}
\begin{picture}(129.60,92.00)
\put(47.00,86.00){\line(1,-1){35.00}}
\put(77.00,56.00){\line(-1,0){30.00}}
\put(67.00,51.00){\line(-1,1){20.00}}
\put(47.00,71.00){\line(1,0){15.00}}
\put(47.00,56.00){\line(1,-1){5.00}}
\put(60.00,6.00){\line(1,0){25.00}}
\put(47.00,21.00){\line(1,0){35.00}}
\put(47.00,21.00){\line(1,-1){15.00}}
\put(77.00,6.00){\line(-1,1){20.00}}
\put(127.00,6.00){\line(-1,0){20.00}}
\put(112.00,6.00){\line(-1,1){5.00}}
\put(127.00,6.00){\line(-1,1){20.00}}
\put(112.00,21.00){\line(-1,0){5.00}}
\put(72.00,26.00){\line(1,-1){10.00}}
\put(47.00,86.00){\circle*{5.20}}
\put(47.00,71.00){\circle*{5.20}}
\put(62.00,71.00){\circle*{5.20}}
\put(62.00,56.00){\circle*{5.20}}
\put(47.00,56.00){\circle*{5.20}}
\put(77.00,56.00){\circle*{5.20}}
\put(77.00,21.00){\circle*{5.20}}
\put(62.00,21.00){\circle*{5.20}}
\put(47.00,21.00){\circle*{5.20}}
\put(62.00,6.00){\circle*{3.20}}
\put(62.00,6.00){\circle{5.20}}
\put(77.00,6.00){\circle*{5.20}}
\put(112.00,6.00){\circle*{5.20}}
\put(127.00,6.00){\circle*{5.20}}
\put(112.00,21.00){\circle*{5.20}}
\put(47.00,82.00){\makebox(0,0)[ct]{$+$}}
\put(47.00,67.00){\makebox(0,0)[ct]{$-$}}
\put(47.00,52.00){\makebox(0,0)[ct]{$+$}}
\put(62.00,67.00){\makebox(0,0)[ct]{$-$}}
\put(62.00,52.00){\makebox(0,0)[ct]{$+$}}
\put(77.00,52.00){\makebox(0,0)[ct]{$+$}}
\put(94.00,21.00){\makebox(0,0)[cc]{$\cdots$}}
\put(94.00,6.00){\makebox(0,0)[cc]{$\cdots$}}
\put(38.00,86.00){\makebox(0,0)[rc]{$c=0$}}
\put(38.00,71.00){\makebox(0,0)[rc]{$c=1$}}
\put(38.00,56.00){\makebox(0,0)[rc]{$c=2$}}
\put(38.00,21.00){\makebox(0,0)[rc]{$c=a-1$}}
\put(38.00,6.00){\makebox(0,0)[rc]{$c=a$}}
\put(62.00,80.00){\makebox(0,0)[cb]{$k=1$}}
\put(77.00,65.00){\makebox(0,0)[cb]{$k=2$}}
\put(112.00,29.00){\makebox(0,0)[cb]{$k=a-1$}}
\put(127.00,15.00){\makebox(0,0)[cb]{$k=a$}}
\put(38.00,40.00){\makebox(0,0)[rc]{$\vdots$}}
\put(94.00,39.00){\makebox(0,0)[cc]{$\ddots$}}
\put(47.00,92.00){\makebox(0,0)[cb]{$k=0$}}
\end{picture}
\caption{\small The extended twisted TPG for $U_q[gl(n|n)^{(2)}]~ (n>2)$
for the tensor product $V(\l_a)
\otimes V(\l_a)$. The vertex labelled by the pair $(c,k)$ corresponds to the
irreducible $U_q[osp(n|n)]$ module $V(k, 2a-2c)$ except for the vertex
corresponding to $c=a,~k=1$, which has been circled to indicate that
it is an indecomposable $U_q[osp(n|n)]$-module.
\label{fig3}}
\end{figure}
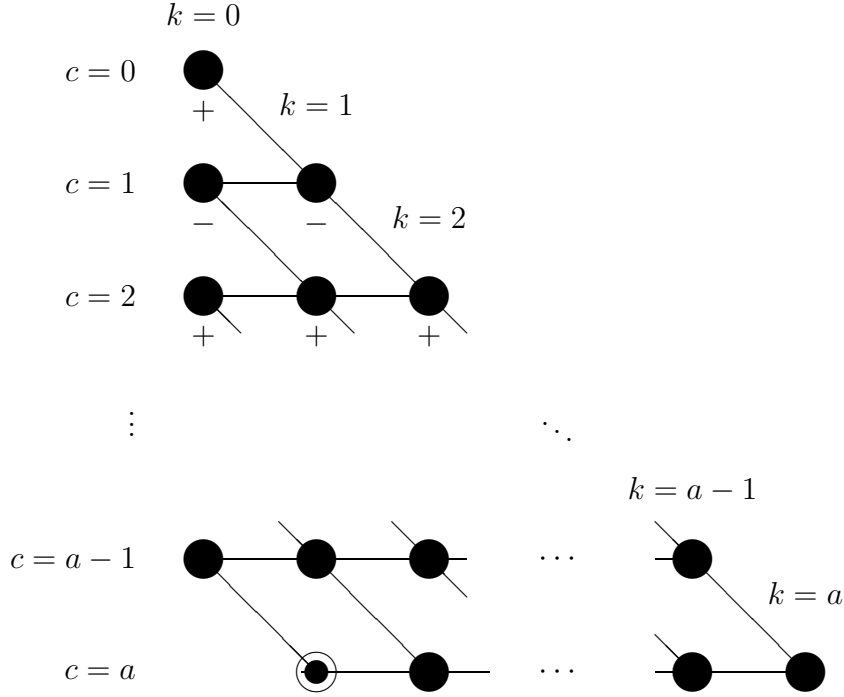

Let $P_V\equiv P^{\l_a\l_a}_V$ be the projection operator from 
$V(\l_a)\otimes V(\l_a)$ onto
$V$ and $P_\nu\equiv P^{\l_a\l_a}_\nu$ the projector onto $V(\nu)$.
Then the R-matrix $\check{R}(z)\equiv \check{R}^{\l_a,\l_a}(z)$
from the extended twisted TPG can be expanded in terms of the operators
$N,~P_V$ and $P_\nu$:
\beq
\check{R}(z)=\rho_N(z) N+\rho_V(z) P_V+\sum_\nu\rho_\nu(z) P_\nu.
    \label{r(z)-check}
\eeq
The coefficients $\rho_\nu(z)$ can be obtained recursively as in the last
subsection from the extended twisted TPG. However, the coefficients $\rho_N(z)$ and 
$\rho_V(z)$ can not be read off from the extended twisted TPG since
the corresponding vertex refers to an indecomposable module. So it
remains to determine these two coefficients. Following our previous
approach \cite{Gou97} to $U_q[gl(2|2)^{(2)}]$, we proceed as follows.

We write the R-matrix $\check{R}\equiv\check{R}(0)=P R,~R\equiv R^{\l_a\l_a}$  
of $U_q[osp(n|n)]$ in the form
\beq
\check{R}=\rho_N(0)N+\rho_V(0)P_V+\sum_\nu\rho_\nu(0) P_\nu,
\eeq
where the coefficients $\rho_N(0),~\rho_V(0),~\rho_\nu(0)$ are all known
from the representation theory of $U_q[osp(n|n)]$,
as is the nilpotent operator satisfying
\beq
P_\nu\, N=N\,  P_\nu=0,~~~~P_V\,N=N\,P_V=N,~~~~N^2=0.\label{N-property}
\eeq
Explicitly, 
\beq
N=\frac{1}{2\rho_V(0)}(R^TR-\rho_V(0)^2)P_V,
\eeq
where $T$ is the usual graded twist map,
\beq
\rho_N(0)=1,~~~~\rho_V(0)=\e_V\,q^{-C_{\l_a}}
\eeq
with $C_{\l_a}$ being the eigenvalue of the universal Casimir of $U_q[osp(n|n)]$
on $V(\l_a)$ and $\e_V$ being the parity of $V$ which is given by $(-1)^a$
(i.e. +1 for $a$ even and -1 for $a$ odd).

Multiplying the Jimbo equation from the right by $P_\nu$ and from the
left by $P_V$, utilising (\ref{r(z)-check}) and (\ref{N-property}), one gets
\bea
&&\lt(\rho_V(z)P_V+\rho_N(z)N\rt)\lt(
z\,e_0\otimes q^{h_0/2}+
q^{-h_0/2}\otimes e_0\rt)P_\nu\no\\
&&
~~~~~=\rho_\nu(z)P_V\left(e_0\otimes q^{h_0/2}+
z\,q^{-h_0/2}\otimes e_0\right) P_\nu.\label{je1}
\eea
Similarly multiplying the Jimbo equation
from the left by $P_\nu$ and from the right by $P_V$
gives
\bea
&&\rho_\nu(z)P_\nu\lt(
z\,e_0\otimes q^{h_0/2}+
q^{-h_0/2}\otimes e_0\rt)P_V\no\\
&&
~~~~~=\rho_V(z)P_\nu\left(e_0\otimes q^{h_0/2}+
z\,q^{-h_0/2}\otimes e_0\right) P_V,\label{je2}
\eea
where we have employed
\beq
P_\nu\lt(e_0\otimes q^{h_0/2}\rt)N=P_\nu\lt(q^{-h_0/2}\otimes e_0\rt)N=0.
   \label{pn}
\eeq
This is seen as follows. We have
$N\lt(V(\l_a)\otimes V(\l_b)\rt)=V(\dot{0}|\dot{0})
\subseteq \bar{V}(\d_1+\d_2)$. Thus 
\beq
\lt(e_0\otimes q^{h_0/2}\rt)N,~~\lt(q^{-h_0/2}\otimes e_0\rt)N~
  \subseteq \bar{V}(\d_1+\d_2)\label{n}
\eeq
from which (\ref{pn}) follows.

Setting $z=0$ into (\ref{je2}) gives
\beq
\rho_\nu(0)P_\nu\lt(
q^{-h_0/2}\otimes e_0\rt)P_V=
\rho_V(0)P_\nu\left(e_0\otimes q^{h_0/2}\rt)P_V.
\eeq
Substituting this equation into (\ref{je2}) we arrive at
\bea
&&\lt(z+\frac{\rho_V(0)}{\rho_\nu(0)}\rt)\rho_\nu(z)P_\nu\lt(
q^{-h_0/2}\otimes e_0\rt)P_V\no\\
&&~~~~~=\lt(1+z\frac{\rho_V(0)}{\rho_\nu(0)}\rt)
\rho_V(z)P_\nu\left(e_0\otimes q^{h_0/2}\rt)P_V.
\eea
Since $P_\nu(e_0\otimes q^{h_0/2})P_V\neq 0$, it follows that
\beq
\rho_V(z)=\frac{z+\rho_V(0)/\rho_\nu(0)}{1+z\,\rho_V(0)/\rho_\nu(0)}\,
  \rho_\nu(z)\label{rho-V}
\eeq
with $\nu$ as in (\ref{nu-node}). Note that since the extended
twisted TPG is consistent, it does not matter which $\nu$ in (\ref{nu-node}) is used.

Before proceeding to the evaluation of $\rho_N(z)$, it is worth noting that
\beq
N\lt(e_0\otimes q^{h_0/2}\rt)N=N\lt(q^{-h_0/2}\otimes e_0\rt)N
  =0.\label{nn}
\eeq
Now multiplying the Jimbo equation from the left and the right by $P_V$
gives rise to
\bea
&&\lt(\rho_N(z)N+\rho_V(z)P_V\rt)
\lt(z\,e_0\otimes q^{h_0/2}+
q^{-h_0/2}\otimes e_0\rt)P_V\no\\
&&
~~~~~=P_V\left(e_0\otimes q^{h_0/2}+
z\,q^{-h_0/2}\otimes e_0\right) \lt(\rho_N(z)N
   +\rho_V(z)P_V\rt),\label{je3}
\eea
which gives, on multiplying from left by $N$ and using (\ref{N-property})
and (\ref{nn}),
\beq
N\lt(z\,e_0\otimes q^{h_0/2}+
q^{-h_0/2}\otimes e_0\rt)P_V
=N\left(e_0\otimes q^{h_0/2}+
z\,q^{-h_0/2}\otimes e_0\right)P_V. 
\eeq
Setting $z=0$ one gets
\beq
N\lt(q^{-h_0/2}\otimes e_0\rt)P_V
=N\left(e_0\otimes q^{h_0/2}\rt)P_V.\label{je4}
\eeq
Multiplying the Jimbo equation from the right by $N$, utilising
(\ref{N-property}) and (\ref{nn}), one has
\beq
P_V\lt(z\,e_0\otimes q^{h_0/2}+
q^{-h_0/2}\otimes e_0\rt)N=
P_V\left(e_0\otimes q^{h_0/2}+
z\,q^{-h_0/2}\otimes e_0\right) N,
\eeq
which leads to, on setting $z=0$,
\beq
P_V\lt(q^{-h_0/2}\otimes e_0\rt)N=
P_V\left(e_0\otimes q^{h_0/2}\rt)N.\label{je5}
\eeq
Also setting $z=0$ into (\ref{je3}) and using (\ref{je4}) and (\ref{je5}),
one gets
\bea
&&\rho_V(0)P_V\lt(e_0\otimes q^{h_0/2}-q^{-h_0/2}\otimes e_0\rt) P_V\no\\
&&~~~~~=N\lt(e_0\otimes q^{h_0/2}\rt)P_V-
        P_V\lt(e_0\otimes q^{h_0/2}\rt)N.\label{je6}
\eea
Finally, substituting (\ref{je4}), (\ref{je5}) and (\ref{je6}) into 
(\ref{je3}) gives
\bea
&&\lt((1+z)\rho_N(z)+(z-1)\frac{\rho_V(z)}{\rho_V(0)}\rt)
    N(e_0\otimes q^{h_0/2})P_V\no\\
&&~~~~~=\lt((1+z)\rho_N(z)+(z-1)\frac{\rho_V(z)}{\rho_V(0)}\rt)
    P_V(e_0\otimes q^{h_0/2})N,\label{je7}
\eea
which can only be satisfied if
\beq
\rho_N(z)=\frac{1-z}{1+z}\cdot\frac{\rho_V(z)}{\rho_V(0)}.\label{rho-N}
\eeq
Indeed, applying (\ref{je7}) to 
$\bar{V}(\d_1+\d_2)\subset V$ we obtain, using $N\bar{V}(\d_1+\d_2)=(0)$,
\beq
\lt((1+z)\rho_N(z)+(z-1)\frac{\rho_V(z)}{\rho_V(0)}\rt)
    N(e_0\otimes q^{h_0/2})\bar{V}(\d_1+\d_2)=(0).
\eeq
Since $N(e_0\otimes q^{h_0/2})\bar{V}(\d_1+\d_2)\neq (0)$, 
one ontains (\ref{rho-N}).

Summarizing, the R-matrix corresponding to (\ref{decom-indecom})
reads
\bea
\check{R}(z)&=&\rho_N(z) N+\rho_V(z)P_V+{\sum_{c=0}^a}'{\sum_{k=0}^c}' 
  \prod^{c-k}_{j=1} \langle 2j-2a \rangle_+\,\no\\
& &  \prod^c_{i=1} \langle i-2a-1 \rangle_-\,
  P_{(2a-2c+k)\d_1+k\d_2},
\eea
where the primes in the sums signify that terms corresponding to
$c=a$ and $k=0,1$ are ommitted from the sums, and $\rho_V(z),~\rho_N(z)$ are given
by 
\bea
\rho_V(z)&=&\frac{z-q^2}{1-zq^2}
   \prod^{a-1}_{j=1} \langle 2j-2a \rangle_+\,
   \prod^{a-1}_{i=1} \langle i-2a-1 \rangle_-,\no\\
\rho_N(z)&=&(-1)^a q^{-a^2}\;\frac{1-z}{1+z}\;\rho_V(z),
\eea
where we have used $\rho_{2\d_1}(0)=(-1)^{a-1}q^{\frac{1}{2} C_{2\d_1}
-C_{\l_a}}$ and $\rho_V(0)=(-1)^aq^{-C_{\l_a}}$.

\sect{Conclusions}

We have shown  how to construct infinite families of new R-matrices with
$U_q[\hat{L}_0=osp(m|n)]$ invariance, arising from the minimal finite
dimensional irreducible representations of the twisted quantum affine
superalgebra $U_q[gl(m|n)^{(2)}]$.

These R-matrices are the only ones so far constructed with $U_q[osp(m|n)]$
invariance apart from the following special exceptions: (i) Those arising
from $\hat{L}_0=osp(2|2)\cong sl(2|1)$, whose R-matrices are already known from
the $gl(m|n)$ case \cite{Del95a}. (ii) Those arising from the vector module of
$U_q(\hat{L}_0)$, which is known to be affinizable to an irreducible
representation of the untwisted quantum affine superalgebra
$U_q[osp(m|n)^{(1)}]$. However the remaining minimal irreducible representations
of $U_q(\hat{L}_0)$ are not affinizable in the untwisted sense, as noted in the
paper. Moreover even in the case of vector representation, the R-matrices
constructed above are different to those arising from the untwisted case.

The R-matrices of this paper will thus give rise to new integrable models
with $U_q[osp(m|n)]$ invariance, which will be investigated elsewhere. It
is particularly interesting in the case $a=b,~m=n>2$, that the R-matrices
admit a $U_q[osp(n|n)]$ invariant nilpotent component, a feature not seen
previously in the untwisted or non-super cases.

\begin{center}
 {\bf Acknowledgements.} 
\end{center}
This paper was completed when YZZ visited Northwest University, China.
He thanks Australian Research Council IREX programme for an Asia-Pacific
Link Award and Institute of Modern Physics of the Northwest University
for hospitality. The financial support from Australian Research 
Council large, small and QEII fellowship grants is also gratefully acknowledged.

\appendix

\sect{Appendix}

Throughout $L_0=o(n=2r)\oplus gl(r)$ denotes the zeroth ${\bf Z}$-graded
component of $\hat{L}_0=osp(m=n|n)$. From \cite{Gou99a}, $V$ in (V.18)
admits a composition series of length 3,
\bea
V\supset \bar{V}(\d_1+\d_2)\supset V(\dot{0}|\dot{0}).\no
\eea
Here $\bar{V}(\d_1+\d_2)$ is indecomposable and cyclically generated by
a $U_q(\hat{L}_0)$ highest weight vector of weight $\d_1+\d_2$ and is the
unique maximal submodule of $V$.
$V(\dot{0}|\dot{0})$ is the trivial one dimensional $U_q(\hat{L}_0)$ module
which is the unique submodule of $\bar{V}(\d_1+\d_2)$, while the factor
module $V/ \bar{V}(\d_1+\d_2)$ is isomorphic, as a $U_q(\hat{L}_0)$ module, to
$V(\dot{0}|\dot{0})$. It is our aim here to prove
\begin{Proposition}\label{prop-A}: The space of $U_q(\hat{L}_0)$-invariant
operator on $V$ is spanned by the identity operator $I$ on $V$ and a
nilpotent invariant operator $N$, unique up to a scalar multiples,
satisfying
\bea
NV=V(\dot{0}|\dot{0}),~~~~N^2=(0).\no
\eea
\end{Proposition}

\noindent{\it Proof.} First we note that $V$ is completely reducible as a
$U_q(L_0)$ module, so we have a $U_q(L_0)$ module decomposition
\beq
V=W\oplus V(\dot{0}|\dot{0})\label{app-1}
\eeq
for some $U_q(L_0)$ submodule $W$ of co-dimension 1 in $V$. Let $v^+_0$ be the
maximal vector of the irreducible $U_q(\hat{L}_0)$ module $V(\dot{0}|\dot{0})$
and $\xi\in V$ the canonical generator of the factor module
$V/ \bar{V}(\d_1+\d_2)$, so that
\beq
(a-\e(a))\xi\in \bar{V}(\d_1+\d_2),~~~~\forall a\in U_q(\hat{L}_0)\label{app-2}
\eeq
with $\e$ being the co-unit.

We define an operator $N$ on $V$ by
\beq
N\xi=v^+_0,~~~~N\bar{V}(\d_1+\d_2)=(0).\label{app-3}
\eeq
Since $\xi$ is uniquely determined modulo $\bar{V}(\d_1+\d_2)$, and $N$ vanishes
on this subspace, $N$ as defined by (\ref{app-3}) is unique, up to scalar
multiples. Moreover $N$ is $U_q(\hat{L}_0)$-invariant since, $\forall
a\in U_q(\hat{L}_0)$, 
\bea
aN\xi=a v^+_0=\e(a)v^+_0=\e(a)N\xi=Na\xi,~~~{\rm by~(\ref{app-2},\;\ref{app-3})},\no
\eea
while it is obvious from (\ref{app-3}) that
\bea
(aN-Na)\bar{V}(\d_1+\d_2)=(0)\Longrightarrow (aN-Na)V=(0),~~~\forall a
  \in U_q(\hat{L}_0),\no
\eea
so $N$ on $V$ is invariant as stated.

Now let $A\in{\rm End}(V)$ be any $U_q(\hat{L}_0)$-invariant and note that
the usual form of Schur's lemma applies to $\bar{V}(\d_1+\d_2)$. Thus there
exists $\a\in{\bf C}$ such that
\bea
(A-\a I)\bar{V}(\d_1+\d_2)=(0),\no
\eea
and in particular, from (\ref{app-2}),
\beq
(A-\a I)(a-\e(a))\xi=(a-\e(a))(A-\a I)\xi=(0),~~~\forall a\in 
   U_q(\hat{L}_0).\label{app-4}
\eeq
In view of the decomposition (\ref{app-1}) we may write
\bea
(A-\a I)\xi=w+\b v^+_0\no
\eea
for some $w\in W,~\b\in{\bf C}$. Hence, $\forall a\in U_q(\hat{L}_0)$,
\bea
\e(a)w+\e(a)\b v^+_0&=&\e(a)(A-\a I)\xi\no\\
&=& a(A-\a I)\xi,~~~{\rm by~(\ref{app-4})}\no\\
&=&aw+\b av^+_0=aw+\e(a)\b v^+_0\no
\eea
$\Longrightarrow aw=\e(a)w,~\forall a\in U_q(\hat{L}_0)$, so
\bea
w\in W\cap V(\dot{0}|\dot{0})=(0).\no
\eea
Hence we must have
\bea
(A-\a I)\xi=\b v^+_0,~~~~(A-\a I)\bar{V}(\d_1+\d_2)=(0)\no
\eea
so that, from (\ref{app-3}),
\bea
A-\a I=\b N\Longrightarrow A=\a I+\b N,\no
\eea
which is sufficient to prove the result.

Finally we note that $N$ acting on $V\subset V(\l_a)\otimes V(\l_a)$
satisfies the requirements of proposition \ref{prop-A}, as desired.


\end{document}